\documentclass[a4paper,10pt]{article}

\usepackage{amsfonts,amssymb,mathrsfs,amscd}

\def \Z {{\mathbb {Z}}}

\def\eps{\varepsilon}

\usepackage[T1]{fontenc} 
\usepackage[cp1251]{inputenc}  
\usepackage[russian]{babel}

\title{\bf Slow convergence almost everywhere of ergodic averages}
\author{\bf Valery V. Ryzhikov}
\date{}
\begin{document}

\maketitle

For an ergodic action of the group $Z^n$ on a probability space and a given arbitrarily slowly decreasing to zero sequence, there exists an integrable function such that the standard ergodic time averages for it converge almost everywhere to the spatial average of the function at a rate that is not asymptotically majorized by this sequence. This generalizes the Krengel effect about the absence of universal estimates for the rate of convergence in Birkhoff's ergodic theorem. The proof uses a weakened version of Rokhlin's lemma for ergodic $Z^n$-actions. It ensures the existence of the required
sequence of asymptotically almost invariant sets with given measures. A feature of the construction of such a sequence is that the choice of the next almost invariant set depends on the original function and on the choice of the previous invariant sets. A significant deviation of  the ergodic averages from the mean of a positive function can be uniformly realized over  extremely large time intervals. The deviation  can be greater than a positive constant and differ little from it on arbitrary long time intervals. We not only can achieve the specified deviations arbitrarily far, but the sequence of such deviations from the average  can be  realized as a sequence wanishing arbitrarily slowly.

\Large
\section{ Медленная сходимость почти всюду эргодических средних}
Пусть  $T$ -- эргодический автоморфизм вероятностного пространства $(X,m)$,   $f\in L_1(X,m)$.  Теорема Биркгофа утверждает, что  для почти всех $x\in X$  выполнено  
$$A(x,N,f):=\frac 1 N \sum_{i=1}^{N} f(T^ix) \to  \int_X f\, dm,\ \ N\to\infty .$$ 
Скорость сходимости   зависит от автоморфизма и функции, в связи с чем  возникло множество разнообразных задач, например, о максимальной скорости сходимости, о связях скоростей с убыванием корреляций, см.  \cite{P}.
Общий результат  был получен Кренгелем \cite{K}:  для  эргодического автоморфизма  $T$ имеет место    сколь угодно медленная скорость сходимости средних  Биркгофа, что  реализуется  при помощи выбора  подходящих функций  $f\in L_1(X,m)$. 
В предлагаемой заметке, модифицируя метод работ  \cite{23}, \cite{25},  мы реализуем  медленную  сходимость средних  Биркгофа для произвольного эргодического $\Z^n$-действия.
 Для простоты ограничесся  рассмотрением  усреднений вдоль временных  (квадратных)  множеств. Обозначим  
$$Q_N=\{z=(z_1,\dots,z_n)\, :\, 1\leq z_1,\dots,z_n\leq N\}.$$
Для эргодического  $\Z^n$-действия $\{T^z\, :\, z\in \Z^n\}$, сохраняющего меру $m$,  и интегрируемой функции $f$ выполняется
$$A(x,N,f):=\frac 1 {|Q_N|} \sum_{z\in Q_N} f(T^zx) \to  \int_X f\, dm,\ \ N\to\infty .$$ 
Это частный случай общей теоремы о поточечной сходимости средних для аменабельных групп \cite{L}.  Покажем, что для заданной сколь угодно медленно  убывающей к 0 последовательности можно как угодно мало изменить по норме  функцию $f$  так, что скорость сходимости средних не будет мажорироваться  этой последовательностью. 

\vspace{2mm} 
\bf Теорема 1. \it Пусть  $\{T^z\}$ -- эргодическое действие группы $\Z^n$ автоморфизмами  вероятностного пространства $(X,m)$ и задан бесконечный набор чисел $a_k>0$ такой, что    $\sum_{i}^\infty a_i<\infty$.  Для всяких неотрицательной $ f_0\in L_1(X,m)$, $\eps>0$ и   последовательности $M_k\to \infty$ найдутся  $N_k>M_k$ и   функция 
$ f\in L_1(X,m)$, $\|f_0-f\|<\eps$  такие, что выполнено 
$$m\left(\ x:\ \left|\, A(x,N_k,f) - \int f \, dm\, \right|\ >\ a_k \ \right)\ \to\ 1.$$   \rm

\vspace{2mm} 
\bf Замечание 1. \rm В \cite{25} рассматривалась похожая задача о медленной сходимости по норме. Для реализации медленной сходимости почти всюду  требуются дополнительные соображения. Забегая вперед, скажем, что для  замедления скорости сходимости выбор очередного почти инвариантного множества $U_{k}$ существенно зависит от выбора предыдущих инвариантных множеств $U_{i}, i<k$.
В случае сходимости по норме (см. \cite{25}) выбор почти инвариантного множества  не зависит от выбора предыдущих.

\vspace{2mm}
Нужные нам почти  инвариантные множества обеспечивает   следствие из  леммы Рохлина  для эргодических $\Z^n$-действий. Последняя является частным случаем общих фактов, установленных Орнстейном и Вейсом \cite{OW}.

\vspace{2mm} 
\bf  Лемма 2. \it  Пусть  $\{T^z\}$ -- эргодическое действие группы $\Z^n$.  Для всяких натурального  $h$ и $a\in (0,1)$ найдется измеримое множество $B$ такое, что  
 $$U=\bigsqcup_{z\in Q_{h}} T^zB_k, \ \ m(U)=a.$$ \rm

 \vspace{2mm} 
Поясним, как  реализуется  замедление сходимости в случае $\Z$-действий. Пусть средние Биркгофа на временных интервалах длины $N_1$  для эргодического автоморфизма и положительной функции $f$ для большинства $x$  стали  близкими к $\int_X f\, dm$.  Обнулив функцию на объединении отрезков траекторий  длины $h\gg N_1$, получим новую функцию $f_1$. В качестве такого объединения используется башня высоты $h$  и  веса $c$ из леммы 2. Тогда вне башни  средние Биркгофа в основном не изменятся, но контролируемо  изменится  интеграл от функции (используем   лемму 3,  которая формулируется ниже).  Таким образом,  мы, так сказать,  замедляем скорость сходимости  для $f_1$   
в момент времени $N_1$. 
Выполнив описанную   процедуру подходящим образом для очень быстро растущей последовательности $N_k\to\infty$,  получаем нужное замедление скорости сходимости временных средних Биркгофа для финальной функции. Последняя по норме может отличаться от исходной функции произвольно мало.
Указанный  метод  работает для  действий  групп $\Z^n$, $n>1$, без принципиальных изменений.

\vspace{2mm} 
\bf Замечание 2. \it Существенное отклонение от средних от интеграла можно равномерно реализовать  на чрезвычайно большом временном интервале.  Например, для большинства $x$  отклонение от среднего может быть больше положительной константы и мало от нее отличаться   для  всех $N\in [N_1, N_1^{N_1^{N_1}}]$.  Для этого надо использовать башни высоты  $h_1\gg N_1^{N_1^{N_1}}$, что разрешается  леммой 2. Таким образом, мы не только достигаем заданных малых 
отклонений как угодно далеко, но и последовательность интегральных отклонений (средняя сумма по всем отклонениям до момента времени $N$) реализуется как сколь угодно  медленно стремящаяся к нулю последовательность. \rm

\vspace{2mm} 
Для контроля значения интеграла от измененной функции нам понадобится следующее утверждение.

\vspace{2mm} 
\bf  Лемма 3. \it Пусть  $\{T^z\}$ -- эргодическое действие группы $\Z^n$ на вероятностном пространстве $(X,m)$.  Если последовательность множеств 
$C_k\subset X$ такова, что $m(C_k)=c>0$  и для вского фиксированного $z$ выполнено $m(T^zC_k\Delta C_k)\to 0$, $k\to\infty$, то для  интегрируемой функции $f$ выполняется  $$\int_{C_k} f\, dm\to c \int_{X} f\, dm, \ \ k\to\infty. $$ \rm

\vspace{2mm}
Утверждение леммы вытекает из классического $\Z^n$-аналога  теоремы фон Неймана  о сходимости средних  $\frac 1 {M^n}\sum_{z\in Q_M} T^{-z} f$   к константе. Действительно, с  учетом асимптотической 
почти инвариантности множеств $C_k$ относительно фиксированных преобразований $T^z$  имеем 
 $$\int_X {\bf 1}_{C_k} f\, dm
 \approx \int_X  \frac 1 {M^n}\sum_{z\in Q_M} {\bf 1}_{T^zC_k}\, f\, dm
  \approx \int_X  {\bf 1}_{C_k}\ \frac 1 {M^n}\sum_{z\in Q_M} T^{-z} f\, dm \approx c \int_X  f\, dm,$$
где \ $\approx $ \ означает, что разность соседних выражений стремится к нулю при   $k\to\infty$ и достаточно медленном стремлении $M$  к бесконечности.

\section{ Доказательство теоремы}
  Пусть заданы $a_i>0$, неотрицательная  функция $f_0\in L_1(X,m)$, причем   выполнено непринципиальное, но технически удобное условие  $$\|f_0\|>1+2\sum_i a_i, \ \ 2\sum_i a_i<1.$$
   Для эргодического $\Z^n$-действия в силу  сходимости почти всюду при $N\to\infty$  средних $A(x,N,f_0)$ к  $\int f_0dm $  найдется число $N_1$ такое, что
$$m\left(x:\ \left|\ A(x,N_1,f_0) - \int f_0dm \ \right|\ <\ a_1/100 \right)\ >1-a_1/100.$$ 
Применяя    леммы    2, 3,   находим   башню 
$U_1$ такую, что  $m(U_1)=2a_1$, причем     для множества $C_1=X\setminus U_1$  выполнено неравенство
$$\left|\ \int  f_0{\bf 1}_{C_1}dm - m(C_1)\int f_0dm \ \right|\ <\ a_1/100.$$
Обозначая $f_1=f_0{\bf 1}_{C_1}$, замечаем, что 
$$\int  f_0dm-\int  f_1dm \ >\ 1.9 a_1\int f_0dm.$$ 
Подходящий выбор достаточно большого числа  $N_1$ и башни $U_1$ достаточной  высоты $h_1\gg N_1$ обеспечивает для всех $x$ вне  определенного ниже  множества $\Delta_1$, $m(\Delta_1)\ll a_1$, равенство
$$A(x,N_1,f_1) =  A(x,N_1,f_0),$$ 
поэтому для таких $x$ выполнено 
 $$\left|A(x,N_1,f_0)- \int f_0\, dm\right| <a_1/10 .$$
Множество $\Delta_1\setminus U_1$ по определению состоит из точек $x$,  для которых $T^zx\in U_1$ для  некоторого
$z\in Q_{N_1}$. Мера этого множества 
не превосходит числа $$m(U_1)\left(\frac {(h_1+N_1)^n} {h_1^n}-1 \right),$$
поэтому  может быть сколь угодно малой при подходящем выборе $h_1$. 

Множество $\Delta_1\cap U_1$ состоит из точек $x\in U_1$,  для которых выполнено $T^zx\notin U_1$ для некоторого  $z\in Q_{N_1}$.  Так как 
$$ m(\Delta_1\cap U_1)=m(U_1)\left( 1- \frac {(h_1- N_1)^n} {h_1^n}\right),$$
 мера множества $\Delta_1\cap U_1$ может быть сделана  сколь угодно малой.
Для   $x\in U_1\setminus\Delta_1$ выполнено 
   $$A(x,N_1,f_1) =  0,$$
 на $U_1\setminus\Delta_1$ отклонение $A(x,N_1,f_1)$ от $\int f_1 dm$ равно $\int f_1 dm> a_1$.
Таким образом, обнулив функцию $f_0$ на множестве $U_1$ меры $2a_1$,  для  $x\notin \Delta_1$ получаем 
отклонение  $A(x,N_1,f_1) - \int f_1> a_1$.
Тем самым  для $k=1$ получено неравенство
$$m\left(x:\ \left|\ A(x,N_k,f_k) - \int f_k\, dm \ \right|\ >\  a_k \right)\ > \ 1-a_k.  \eqno  (k)$$
Совершенно аналогично  при $k=2, 3,\dots$ для  числа $a_k>0$  и    функции $f_{k-1}$,   $\|f_{k-1}\|>1+ 2\sum_{i=k} a_i$, мы поочередно находим   $N_k$, определяем новую функцию $f_k$, обнуляя $f_{k-1}$ на подходящей  башне $U_k$. Быстрый рост $N_k\to \infty$ и   $h_k/N_k\to \infty $  влечет    выполнение условий $(k)$.

Для финальной функции 
$$ f=f_0{\bf 1}_C, \ \ C=X\setminus \bigcup_{k= 1}^{\infty} U_{k},$$  при достаточно  быстом росте размеров $h_k$ башен $U_k$ меры множеств $\Delta_k$ (определяются аналогично множеству $\Delta_1$)  очень малы, поэтому  выполняются  неравенства 
$$m\left(x:\ \left|A(x,N_k,f) - \int f\, dm \right|>\  a_k \ \right)\ > \
1-2\sum_{i=k}^\infty a_i.$$

Последовательность $N_k$ не ограничена сверху, поэтому  
для  неотрицательной $f\neq 0$ мы получили  для почти всех $x$ требуемую   медленную скорость сходимости $A(x,N_k,f)$ к  $\int f\, dm$. \rm 
Техническое ограничение   $\|f_0\|>1+2\sum_i a_i$, $2\sum_i a_i<1$ не играет 
принципиальной роли, общий случай из него вытекает.  В утверждении теоремы сказано, что функция $f$ может как угодно мало отличаеться по норме от заданной $f_0$. Это достигается, например,  изменением  последовательности $a_i$ на конечном множестве индексов $i$ таким, что $2\sum_{i=k}^\infty a_i<\eps.$   Тогда соответствующая последовательности  $a_i$ финальная функция  $f$ удовлетворяет условию  $\|f_k-f\|<\eps$ и  для нее  имеет место требуемая  медленная скорость сходимости средних. Теорема доказана.

\vspace{2mm}
\bf Замечание 3. \rm Идею обнуления неотрицательной функции на последовательности 
асимптотически почти инвариантных множеств можно использовать  для 
решения аналогичной задачи в случае действий счетных  аменабельных групп.  
Для сходимостей по норме эта идея реализована  в \cite{25}.   Эффект медленной поточечной сходимости для действий аменабельных групп будет предствлен  в отдельной заметке И.В. Бычкова.

\end{document}